\input epsf
\documentclass[a4paper,11pt]{amsart}
\title{An Orchard Theorem}
\author{Roland Bacher
}
\date{}

\begin{document}
\maketitle

{\begin{abstract} We describe a natural way to plant
cherry- and plumtrees at prescribed generic locations
in an {orchard.}
\end{abstract} }\footnotetext{Support from the Swiss National 
Science Foundation is gratefully acknowledged.}


\section{Main results}

The main result of this paper may be paraphrased comprehensively as
follows: Most people would agree that a natural way to plant trees
of two species along a row is to alternate them.
Our main result (the Orchard Theorem) generalises this to higher
dimensions. In dimension 2 it implies that there is a natural way to
plant cherrytrees and plumtrees
in an orchard if the prescribed
locations of the trees are generic (no alignments of 
three trees). Figure 1 shows such an orchard planted
with 3 cherry- and 6 plumtrees (they play of course symmetric
roles).
\medskip

\centerline{\epsfysize8cm\epsfbox{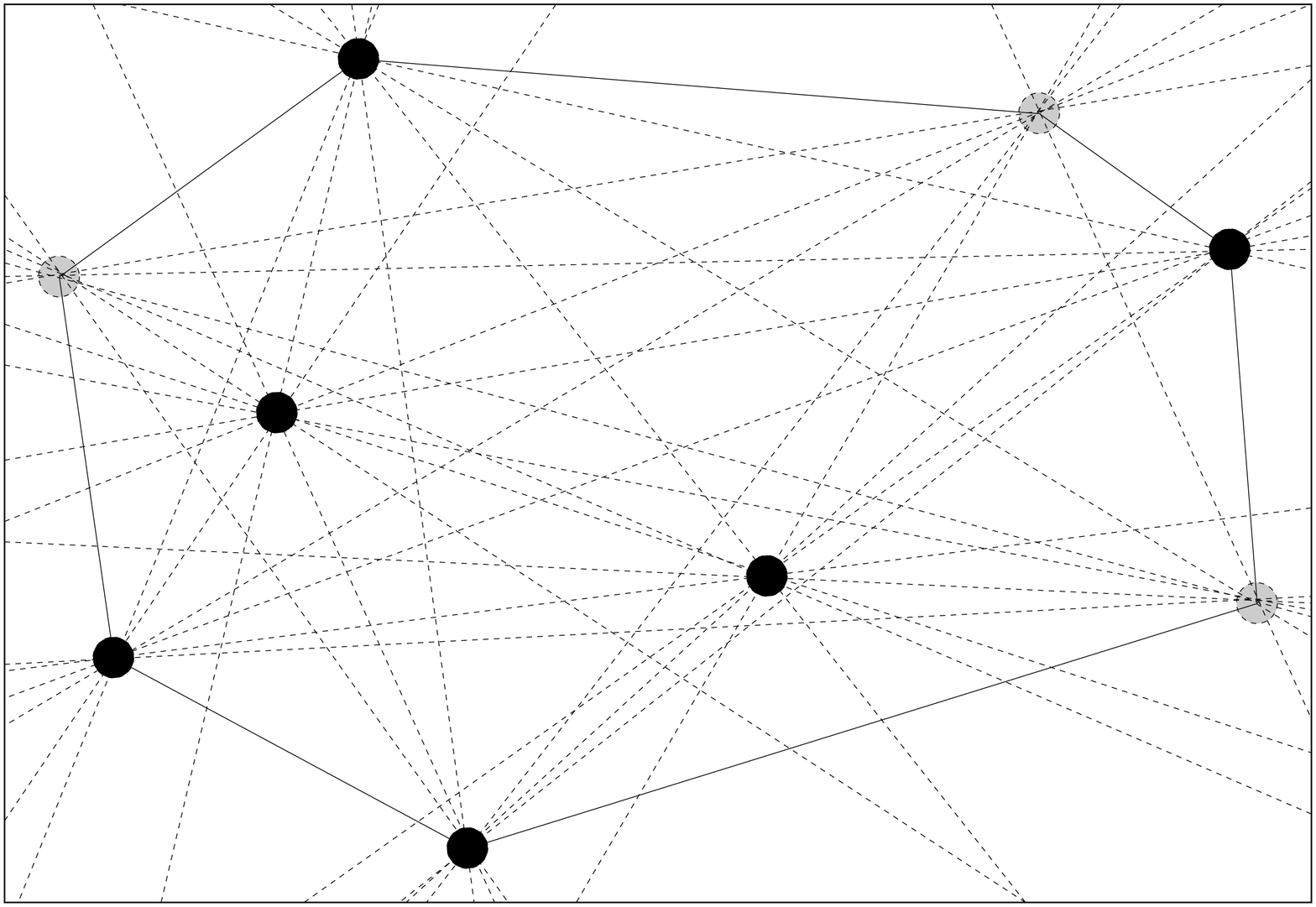}}
\centerline{Figure 1:
An orchard having 3 cherry- and 6 plumtrees in generic positions} 
\medskip

A finite set ${\mathcal P}=\{P_1,\dots,P_n\}$ of $n$ points in the 
oriented real affine space
${\bf R}^d$ is a {\it generic configuration} if any subset of $k+1\leq d+1$
points is affinely independent. Generic configurations of $n\leq d+1$
points in ${\bf R}^d$ are simply vertices of an $(n-1)-$dimensional
simplex. For $n\geq d+1$, 
genericity boils down to the fact that any set of $d+1$ 
points in ${\mathcal P}$ spans ${\bf R}^d$ affinely. 

Two generic configurations ${\mathcal P}^1$ and ${\mathcal P}^2$ are 
{\it isomorphic} if there
exists a bijection $\varphi:{\mathcal P}^1\longrightarrow {\mathcal P}^2$
such that all corresponding $d-$dimensional simplices (with vertices 
$(P_{i_0},\dots,
P_{i_d})\subset {\mathcal P}^1$ and $(\varphi(P_{i_0}),\dots,
\varphi(P_{i_d}))\subset {\mathcal P}^2$) have the same orientation
(given for instance for the first simplex 
by the sign of the determinant of the $d\times d$ matrix with rows 
$P_{i_1}-P_{i_0},\dots,P_{i_d}-P_{i_0}$).

Two generic configurations ${\mathcal P}(-1)$ and ${\mathcal P}(+1)$ 
are {\it isotopic} if there exists a continuous path (with respect to
the obvious topology on ${\bf C}^{dn}$)
of generic configurations $t\longmapsto {\mathcal P}(t)$
which joins them. Isotopic
configurations are of course isomorphic. I ignore to what extend the
converse holds. 

An affine hyperplane $H\subset {\bf R}^d$ {\it separates} two points $P,Q\in
{\bf R}^d\setminus H$ if $P,Q$ are not in the same connected component of 
${\bf R}^d\setminus H$.
For two points $P,Q$ of a generic configuration 
${\mathcal P}=\{P_1,\dots,
P_n\}\subset {\bf R}^d$ we denote by $n(P,Q)$ the number of distinct 
hyperplanes separating $P$ and $Q$ which are
affinely spanned by $d$ distinct elements 
in ${\mathcal P}\setminus \{P,Q\}$. The
number $n(P,Q)$ depends obviously only of the isomorphism type of 
$\mathcal P$.

{\bf Theorem 1.1 (Orchard Theorem).} {\sl The relation defined by $P\sim Q$ if either
$P=Q$ or if $$n(P,Q)\equiv  {n-3\choose d-1}\pmod 2$$ is an
equivalence relation having at most 2 classes 
on a generic configuration ${\mathcal P}$ of $n$ points in ${\bf R}^d$.}

We call the equivalence relation of Theorem 1.1 the {\it Orchard
relation} and the induced partition on $\mathcal P$ the
{\it Orchard partition}.

{\bf Example 1.2.} Consider a configuration 
${\mathcal P}\subset {\bf S}^2\subset {\bf R}^3$ consisting of $n$ points 
contained in the Euclideean unit sphere ${\bf S}^2$ and which are generic
as a a subset of ${\bf R}^3$ in the above sense, i.e. $4$ distinct 
points of ${\mathcal P}$ are never contained in a common affine plane
of ${\bf R}^3$. A stereographic projection $\pi: {\bf S}^3\setminus \{N\}
\longrightarrow {\bf R}^2$
with respect to a point $N\in {\bf S}^2\setminus {\mathcal P}$
sends the set ${\mathcal P}\subset {\bf S}^2$ into a set
$\tilde{\mathcal P}=\pi({\mathcal P})\subset {\bf R}^2$ such that $4$ points
of $\tilde {\mathcal P}$ are never contained in a common Euclideean
circle or line of ${\bf R}^2$. The Orchard relation on ${\mathcal P}$ 
can now be seen on $\tilde {\mathcal P}$ as follows: Given two distinct 
points $\tilde P\not=\tilde Q\in\tilde{\mathcal P}$ count the number
$n(\tilde P,\tilde Q)$ of circles or lines determined by $3$ points
in $\tilde {\mathcal P}\setminus\{\tilde P,\tilde Q\}$ which separate them.
The points $P,Q\in{\mathcal P}$ are now Orchard-equivalent if and only if
$n(\tilde P,\tilde Q)\equiv {n-3\choose 2}\pmod 2$.

A consequence of the Orchard Theorem is the fact that 
points of a generic configuration ${\mathcal P}
\subset {\bf R}^d$ carry a structure which can be encoded by a
rooted binary tree $\mathcal T$: 
Vertices of $\mathcal T$ are suitable subsets of $\mathcal P$ and define 
hence generic subconfigurations of ${\mathcal P}$. 
The root corresponds to the complete set $\mathcal P$. The remaining 
vertices of ${\mathcal T}$ are defined recursively as follows:
The two sons (if they exist) of a vertex $V_{\tilde{\mathcal P}}\in 
{\mathcal T}$
corresponding to a subset ${\tilde{\mathcal P}}\subset {\mathcal P}$
are the two non-empty Orchard classes (equivalence classes) 
of the generic configuration $\tilde{\mathcal P}$.

It would of course be interesting to understand the leaves (or atoms)
of such trees. They correspond to configurations consisting only of
equivalent points. Generic configurations of ${\bf
  R}^d$ having at most $d+1$ elements are of course such leaves but
there are many others (e.g. vertices of a convex plane polygone having
an odd number of vertices). 

A {\it flip} is a continuous path $$t\longmapsto
{\mathcal P}(t)=(P_1(t),\dots,P_n(t))\in\left({\bf R}^d\right)^n,\ 
t\in [-1,1]$$
with ${\mathcal P}(t)=\{P_1(t),\dots,P_n(t)\}$ generic except for $t=0$
where there exists exactly one set ${\mathcal F}(0)=
(P_{i_0}(0),\dots,P_{i_d}(0))$, called the {\it flipset},
of $(d+1)$
points contained in an affine hyperplane spanned by any subset of
$d$ points in
${\mathcal F}(0)$. We require moreover that the simplices
$(P_{i_0}(-1),\dots,P_{i_d}(-1))$ and $(P_{i_0}(1),\dots,P_{i_d}(1))$
carry opposite orientations. Geometrically this means that a point 
$P_{i_j}(t)$ crosses the hyperplane spanned by
${\mathcal F}(t)\setminus \{P_{i_j}(t)\}$ at time $t=0$.

It is easy to see that two generic configurations ${\mathcal
P}^1,{\mathcal P}^2\subset {\bf R}^d$ having $n$ points
can be related by a continuous path involving at most a finite 
number of flips.

The next result shows that flips modify the Orchard relation
only locally.

{\bf Proposition 1.3 (Flip Proposition).} {\sl Let ${\mathcal P}(-1),{\mathcal P}(+1)\subset 
{\bf R }^d$ be two generic 
configurations related by a flip with respect to a subset 
${\mathcal F}(t)$ of $(d+1)$ points. 

\ \ (i) If two distinct points $P(t),Q(t)$ are either both contained in 
${\mathcal F}(t)$ or both contained in its complement
${\mathcal P}(t)\setminus {\mathcal F}(t)$ then we have 
$$P(-1)\sim Q(-1)\hbox{ if and only if }P(+1)\sim Q(+1)\ .$$

\ \ (ii) For $P(t)\in {\mathcal F}(t)$ and $Q(t)\not\in {\mathcal
  F}(t)$ we have
 $$P(-1)\sim Q(-1)\hbox{ if and only if }P(+1)\not\sim Q(+1)\ .$$}

\medskip

\centerline{\epsfysize6cm\epsfbox{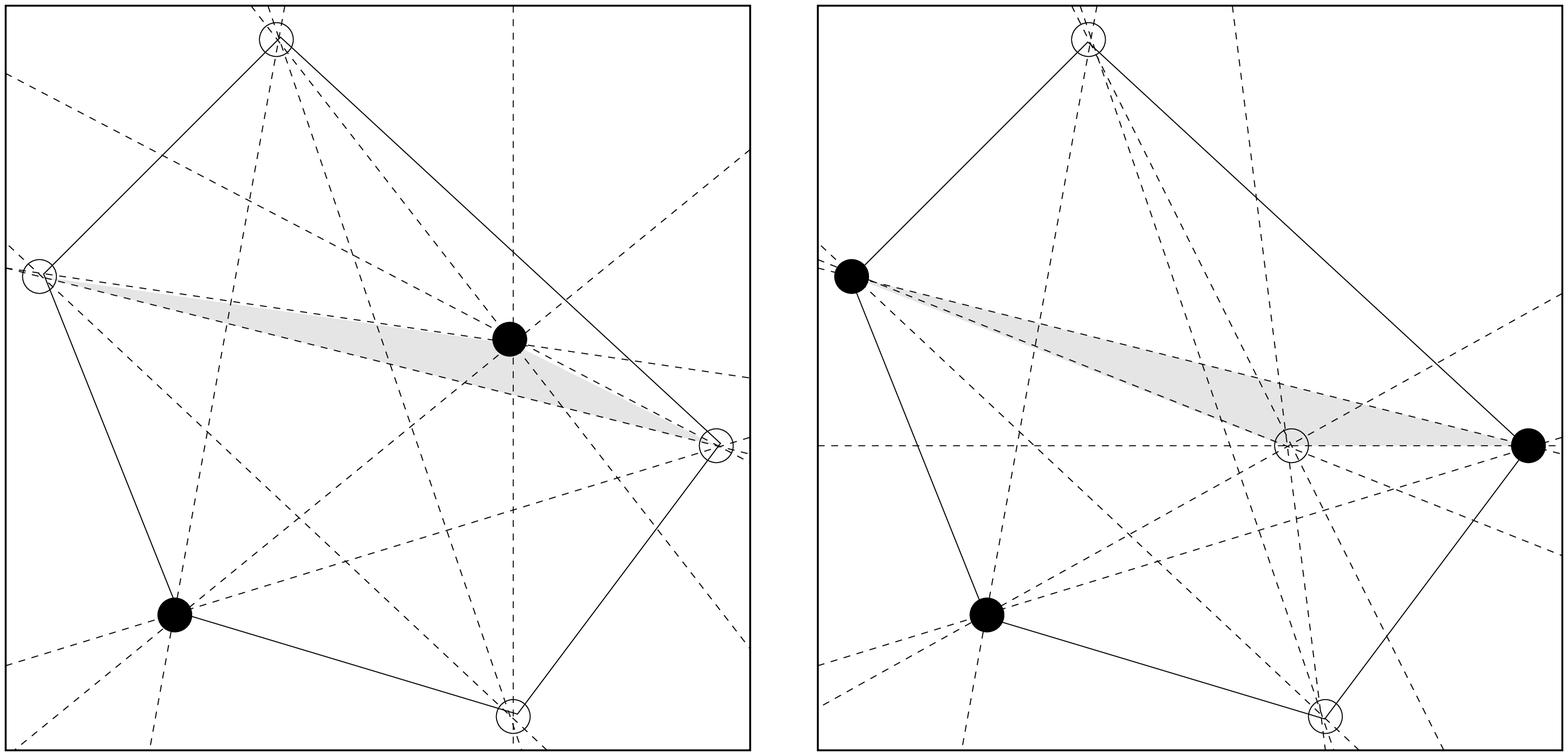}}
\centerline{Figure 2:
Two orchards with 6 trees related by a flip} 
\medskip

The Flip Proposition bears bad news for a possible natural
generalisation of the Orchard relation to non-generic configurations. 
The two equivalence classes  play indeed a totally symmetric role with
respect to the flipset and there seems no natural
way to break this symmetry for the non-generic configuration
${\mathcal P}(0)$ involved in a flip. It is however possible
to define the Orchard relation on {\it generic points} of 
arbitrary configurations: Call a point $P\in {\mathcal P}\subset {\bf R}^d$
of a subset of points {\it generic} if the affine span of $\{P,P_{i_1},\dots,
P_{i_k}\}$ is $k-$dimensional (for $k\leq d$)
for all subsets $\{P_{i_1},\dots,P_{i_k}\}\subset
{\mathcal P}\setminus \{ P\}$. A small generic perturbation 
$\tilde {\mathcal P}$ of $\mathcal P$ allows then to compute the
orchard relation on the set of generic points of $\mathcal P$ by
considering the restriction of the Orchard relation on $\tilde{
\mathcal P}$ to the image of generic points.

However, the Flip Proposition suggests also perhaps interesting
problems concerning generic configurations: 
Call two generic
configurations of $n$ points in ${\bf R}^{2d+1}$ {\it
orchard-equivalent} if they can be related by a series of flips
whose flipsets have always exactly $(d+1)$ points in each class.
  
More generally, flips are of different types according to
the number of points of each class involved in the corresponding
flipset. A very special type of flips are the {\it
  monochromatic} ones, defined as involving only vertices of one
class in their flipset. 

Understanding isotopy (or more generally isomorphism) classes
of generic configurations
up to flips subject to some restrictions (e.g. only
monochromatic flips or configurations up to
orchard-equivalence in odd dimensions)
might be an interesting problem.

Before describing a last consequence of the Flip Proposition we need
a few notations:

Let ${\mathcal P}$ be a generic configuration with equivalence classes
${\mathcal A}({\mathcal P})$ and ${\mathcal B}({\mathcal P})$
consisting of $a({\mathcal P})$ and $b({\mathcal P})$ elements. 
Call a generic
configuration {\it pointed} if one of the equivalence classes, say
${\mathcal A}({\mathcal P})$ has been selected and denote by 
$\overline{a}({\mathcal P})$ the number of elements in the selected 
class. Let us moreover introduce the finite graph
$\Gamma^p{\mathcal C}_n({\bf R}^d)$ with vertices isomorphism classes
of pointed generic configurations of $n$ points in ${\bf R}^d$, 
two edges beeing joined by an edge if the corresponding
pointed configurations can be related by a flip (a flip of a pointed
cofigurations does not change the selected equivalence class outside
the flipset). The following result is then an easy consequence of the
Flip Proposition.

{\bf Corollary 1.4.} {\sl (i) If $d$ is even then the graph
$\Gamma^p{\mathcal C}_n({\bf R}^d)$ is bipartite, the class of a
pointed configuration $\mathcal P$ beeing given by ${\overline a}(
{\mathcal P})\pmod 2$.

\ \ (ii) If $d$ is odd and $n=2m+1$ is odd then the 
graph $\Gamma^p{\mathcal C}_n({\bf R}^d)$ has two connected
components, the function ${\overline a}({\mathcal P})\pmod 2$
beeing constant on each component.

\ \ (iii) If $d$ is odd and $n=2m$ is even 
then there exist $\pi(n,d)\in \{0,1\}$ such
that ${\overline a}({\mathcal P})\equiv {\overline b}
({\mathcal P})\equiv \pi(n,d)\pmod 2$
for every generic configuration $\mathcal P$ of $n=2m$ points in ${\bf
  R}^d$.}

In the case $d=1$, it is easy to see that the function $\pi(2m,1)$ of
assertion (iii) is given by $\pi(2m,1)\equiv m\pmod 2$. More generally
(cf. Example 3.2 of Section 3), one has
$$\pi(2m,2d+1)\equiv \left\{\begin{array}{ll}
m\pmod 2&\hbox{if }{2m-3\choose 2d}\hbox{ is odd}\cr
0&\hbox{otherwise.}\end{array}\right.$$

Call an antipodal subset ${\mathcal P}=\{\pm P_1,\dots,\pm P_n\}\subset
{\bf S}^{d}$ of the $d-$dimensional unit sphere ${\bf S}^d\subset {\bf R}^{d+1}$
{\it generic} if the linear span in ${\bf R}^{d+1}$ of any subset $\pm P_{i_1},
\dots,\pm P_{i_k}$ of $k$ pairs of points is of dimension $k$ for $k\leq d+1$.

{\bf Theorem 1.5 (Spherical Orchard Theorem).} {\sl 
There exists a natural equivalence relation
having at most $2$ classes on the set ${\mathcal P}=\{\pm P_1,\dots,\pm P_n\}
\subset {\bf S}^d$ of points of a generic antipodal spherical configuration.

If ${n-2\choose d}$ is even, this relation satisfies
$$-P_i\sim P_i\hbox{ for }P_i\in {\mathcal P}$$
and if ${n-2\choose d}$ is odd we have
$$-P_i\not\sim P_i \hbox{ for }P_i\in {\mathcal P}\ .$$}

We have of course also the obvious version (involving suitable
subsets of $d+1$ pairs of antipodal points) of the Flip proposition.

Let ${\mathcal C}$ be a set of continuous real functions on ${\bf R}^k$. 
Suppose ${\mathcal C}$ is a $(d+1)-$dimensional vector
space containing the constant functions. We denote by ${\mathcal C}_0$ the
$d-$dimensional subspace 
$${\mathcal C}_0=\{f\in{\mathcal C}\ \vert\ f(0)=0\}\ .$$
Call a set ${\mathcal P}\subset {\bf R}^k$ of $n$ points 
{\it ${\mathcal C}-$generic} if for each subset $S=\{P_{i_1},\dots,
P_{i_d})$ of $d$ distinct points in $\mathcal P$ if the set
$$I(S)=\{f\in{\mathcal C}\ \vert\ f(P_{i_j})=0,\ j=0\dots,d\}$$
is $1-$dimensional all ${n\choose d}$ lines in ${\mathcal C}$ of this
form are distinct.

Given $P,Q\in{\mathcal P}$, call a set $S=\{P_{i_1},\dots,
P_{i_d})$ of $d$ points as above {\it ${\mathcal C}-$separating}
(or {\it separating} for short) if $f(P)f(Q)<0$ for any $0\not=f\in I(S)$
and denote by $n_{\mathcal C}(P,Q)$ the number of ${\mathcal C}-$separating
subsets of $\mathcal P$.

{\bf Corollary 1.6.} {\sl The relation $P\sim_{\mathcal C} Q$ if
either $P=Q$ or
$$n_{\mathcal C}(P,Q)\equiv {n-3\choose d-1}\pmod 2$$
defines an equivalence relation having at most two classes on a 
set ${\mathcal P}=\{P_1,\dots,P_n\}\subset {\bf R}^k$ 
of $n$ points in ${\bf R}^k$ which are ${\mathcal C}-$generic.}

{\bf Examples 1.7.} (i) Considering the $(d+1)-$dimensional 
vector space of all affine functions in ${\bf R}^d$, Corollary 1.6
boils down to Theorem 1.1.

\ \ (ii) Consider the vector space $\mathcal C$ of all polynomial
functions ${\bf R}^2\longrightarrow {\bf R}$. A finite subset $
{\mathcal P}\subset {\bf R}^2$ is ${\mathcal C}-$generic if and only if
every subset of five points in $\mathcal P$ defines a unique conic and
all these conics are distinct.

\ \ (iii) Consider the vector space $\mathcal C$ of all polynomials of degree
$< d$ in $x$ together with the polynomials $\lambda y,\ 
\lambda\in{\bf R}$. A subset ${\mathcal P}=\{(x_1,y_1),\dots,(x_n,y_n)\}$ 
with $x_1<x_2<\dots,<x_n$ is ${\mathcal C}-$generic if all ${n\choose 
d}$ interpolation polynomials in $x$ defined by $d$ points of
$\mathcal P$ are distinct.

The sequel of this paper is organized as follows:

The next section contains the (elementary) proofs of Theorem 1.1, 
Proposition 1.3 and Corollaries 1.4, 1.6.

Section 3 contains a brief description for concrete computations of
the Orchard relation.

Section 4 states the projective version of the Orchard Theorem and contains the 
proof (and construction of $\sim$) for Theorem 1.5.

Section 5 is devoted to simple arrangements of pseudolines in ${\bf R}P^2$.

\section{Proofs}

In this section we give a proof of Theorem 1.1, Proposition 1.3 and 
Corollary 1.4.

The case $d=1$ (of both results) is of course trivial: A generic 
configuration of $n$ points in ${\bf R}$ is simply a strictly
increasing sequence $P_1<P_2<\dots<P_n$ of $n$ real numbers
and the Orchard Theorem restates the obvious fact that we get an
equivalence relation by considering
$P_i\sim P_j$ for $i\equiv j\pmod 2$. 

We suppose now $d\geq 2$ and introduce some useful notations.

Consider $3$ points $P_i,P_j,P_k$ of a generic configuration
${\mathcal P}=\{P_1,\dots,P_n\}\subset {\bf R}^d$.
If $d\geq 2$ which we will suppose in the sequel, the affine span of
the points $P_i,P_j,P_k$ is a $2-$dimensional affine plane. 
The three affine lines $L_{s,t}$ spanned by two
points $P_s\not=P_t\in \{P_i,P_j,P_k\}$ contain three compact intervals 
denoted $[P_s,P_t]$ and
subdivide the projective plane $\overline\Pi\supset \Pi$
into four triangles $\Delta_0,
\Delta_i,\Delta_j,\Delta_k$ as shown in Figure 3.
\medskip

\centerline{\epsfysize4cm\epsfbox{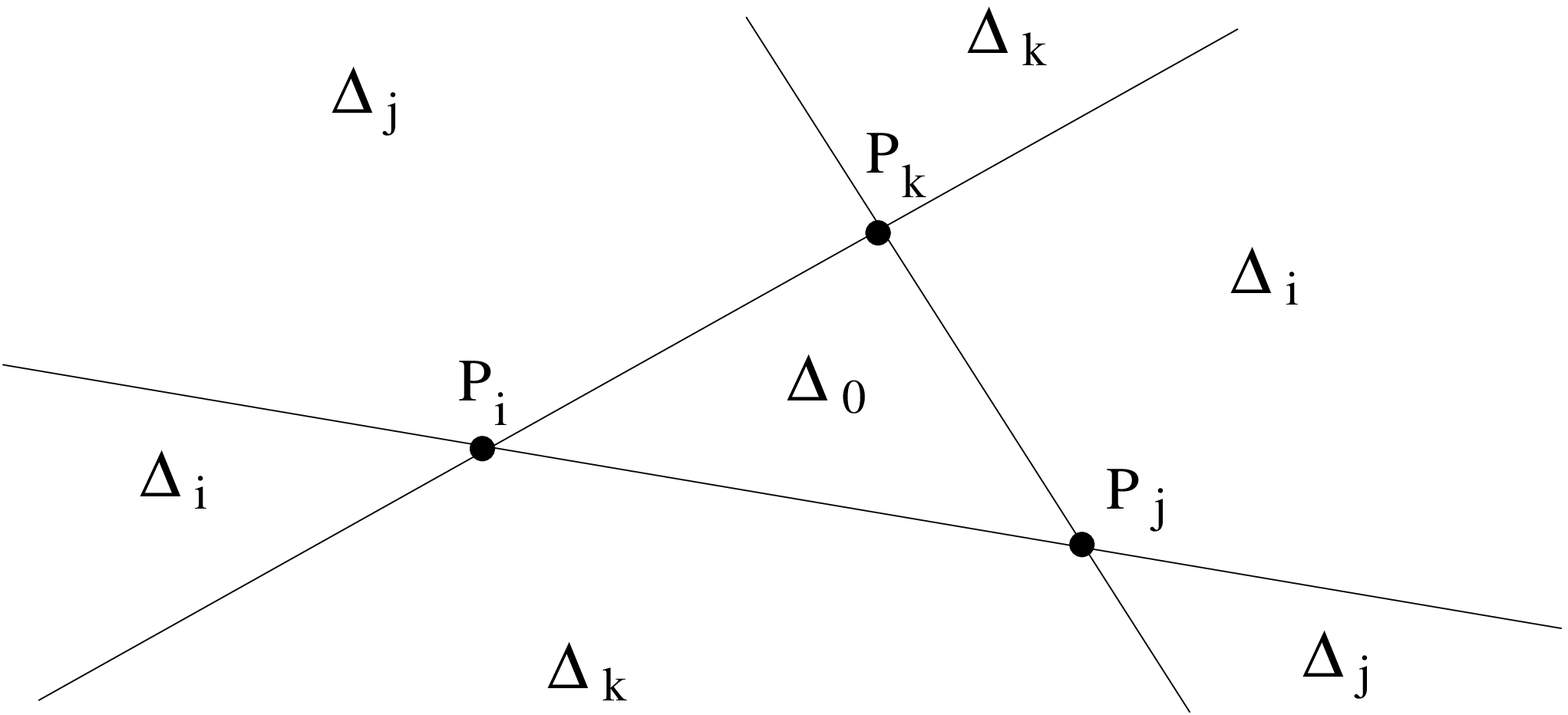}}
\centerline{Figure 3:
The projective plane $\overline{\Pi}=\Delta_0\cup\Delta_i\cup\Delta_j\cup\Delta_k$} 
\medskip
Let $\alpha_i$ denote the number of hyperplans containing
$d$ distinct points in 
${\mathcal P}\setminus \{P_i,P_j,P_k\}$ which intersect both
segments $[P_i,P_j]$ and $[P_i,P_k]$. Introduce $\alpha_j$ and
$\alpha_k$ similarly.

For a subset 
$$S=\{P_{i_1},\dots,P_{i_{d-1}}\}\subset {\mathcal P}
\setminus \{P_i,P_j,P_k\}$$
of $(d-1)$ distinct points in ${\mathcal P}
\setminus \{P_i,P_j,P_k\}$, we denote by $P_S$ the intersection
(which might be at infinity) of the projective span of $S$ with the 
projective plane $\overline \Pi$. We claim that $P_S$ consists of
exactly one point contained in the interior of exactly one of the
four triangles $\Delta_0,\Delta_i,\Delta_j,\Delta_k$: Indeed,
$P_S$ is non-empty by a dimension argument. If $P_S$ is of
dimension $>0$ or if $P_S$ is included in the projective line
$\overline{L_{i,j}}$, the set $S\cup\{P_i,P_j\}$ of $(d+1)$ distinct 
points in $\mathcal P$ is affinely dependent. The same argument holds
of course for the projective lines $\overline{L_{i,k}}$ and 
$\overline{L_{j,k}}$.
 
For $*\in \{0,i,j,k\}$,  we can hence introduce the number 
$\sigma_*$ of subsets 
$$S=\{P_{i_1},\dots,P_{i_{d-1}}\}\subset {\mathcal P}
\setminus \{P_i,P_j,P_k\}$$
of $(d-1)$ distinct points in ${\mathcal P}\setminus \{P_i,P_j,P_k\}$
whose projective span intersects the projective plane $\overline\Pi$ 
in an interior point $P_S$ of the triangle $\Delta_*$.

Recall that $n(P_i,P_j)$ (and similarly $n(P_j,P_k),n(P_i,P_k)$)
denotes the number of hyperplanes spanned by $d$ points in 
${\mathcal P}\setminus\{P_i,P_j\}$ which separate $P_i$ from $P_j$
in the affine space ${\bf R}^d$.

{\bf Lemma 2.1.} {\sl We have
$$\begin{array}{l}
\displaystyle n(P_i,P_j)=\alpha_i+\alpha_j+\sigma_0+\sigma_k\ ,\cr
\displaystyle n(P_j,P_k)=\alpha_j+\alpha_k+\sigma_0+\sigma_i\ ,\cr
\displaystyle n(P_i,P_k)=\alpha_i+\alpha_k+\sigma_0+\sigma_j\ .\end{array}$$}

{\bf Proof of Lemma 2.1.} We prove the formula for $n(P_i,P_j)$. The
remaining cases follow by symmetry.

A hyperplane $H$ separating $P_i$ from $P_j$ 
intersects the interior of the segment $[P_i,P_j]$ and we have 
two subcases depending on the position of $P_k$ with respect to $H$.

If $P_k\not\in H$, the line $H\cap \Pi$ cuts the interior of
the edge $[P_i,P_j]\subset \Delta_0$ and leaves the triangle
$\Delta_0$ by crossing the interior of either the edge $[P_i,P_k]$ or of the
edge $[P_j,P_k]$. Such a hyperplane $H$ contributes hence $1$ to either 
$\alpha_i$ or $\alpha_j$. The line $H_3$ of Figure 3
shows the intersection of such a hyperplane with the plane 
$\Pi$. It yields a contribution of $1$ to $\alpha_j$.

In the remaining case $P_k\in H$, the projective hyperplane
$H$ is spanned by $P_k$ and by a subset $S$ consisting of exactly 
$(d-1)$ points in ${\mathcal P}
\setminus \{P_i,P_j,P_k\}$. The projective line $\overline{\Pi\cap H}$ is hence
defined by the point $P_k\in \Pi$ and by the point $P_S$ (which might
very well be located at infinity), defined
as above as the intersection of the projective span of $S$ with the
projective plane $\overline\Pi$. Since the line $\Pi\cap H$ crosses
$[P_i,P_j]$, the point $P_S\in \overline{\Pi}$ 
belongs either to the interior of
$\Delta_0$ or $\Delta_k$ (it cannot be on the boundary by genericity,
see above) and such a hyperplane $H$ yields hence a
contribution of $1$ to either $\sigma_0$ or $\sigma_k$. 
Figure 4 shows two such situations giving rise to hyperplanes intersecting
$\Pi$ along the lines $H_1$ and $H_2$. The corresponding points 
$P_{S_1}$ and $P_{S_2}$ belong respectively to $\Delta_k$ and
$\Delta_0$.
\hfill QED

\medskip

\centerline{\epsfysize6cm\epsfbox{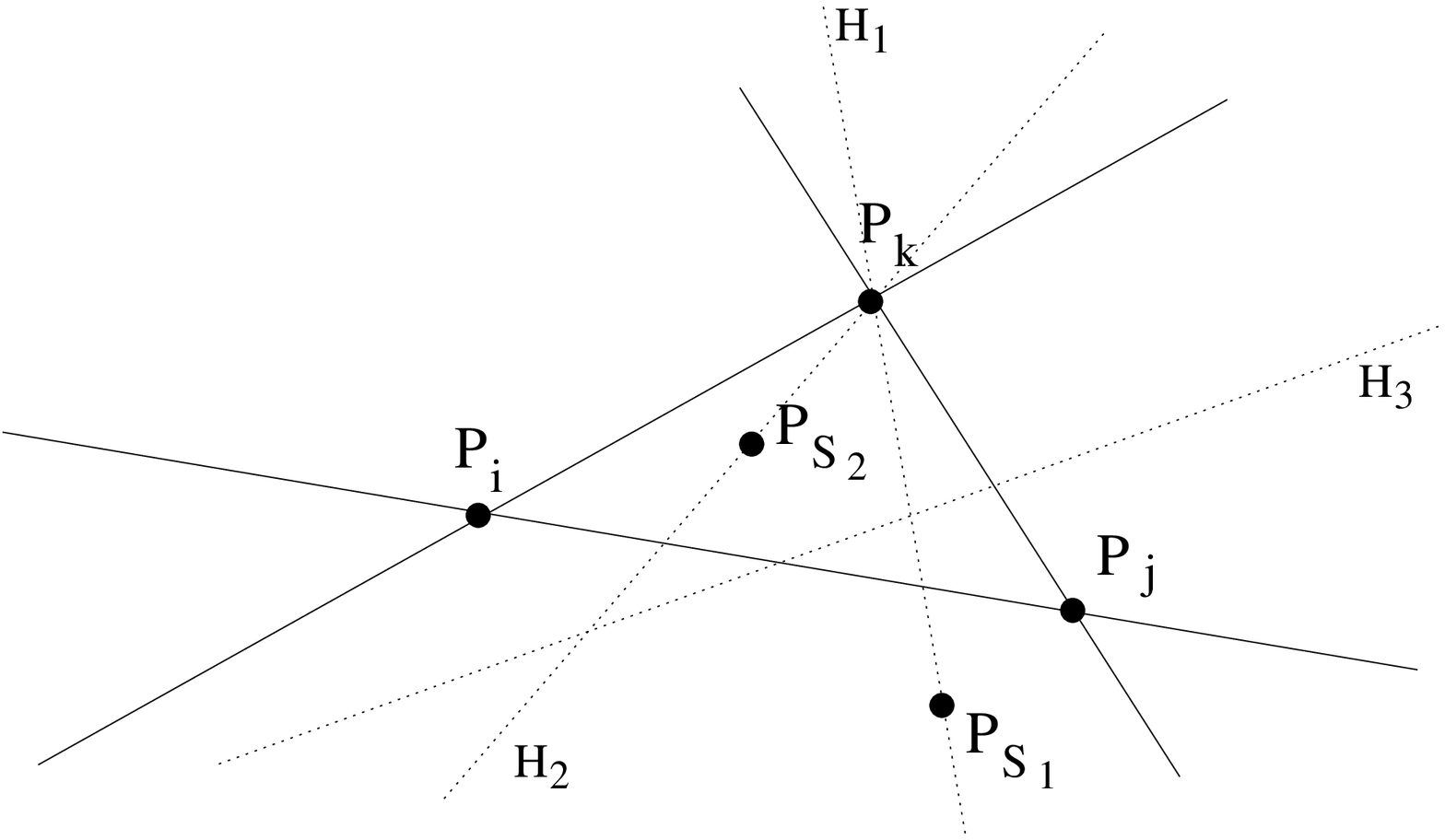}}
\centerline{Figure 4: Three hyperplanes contributing to $n(P_i,P_j)$} 
\medskip

{\bf Proof of Theorem 1.1 in the case $d\geq 2$.} It is obvious that
the Orchard relation $\sim$ is reflexive and symmetric. Transitivity is 
however not completely obvious. We consider hence three points $P_i,P_j$ and 
$P_k$ satisfying $P_i\sim P_j$ and $P_j\sim P_k$.

We have then by Lemma 2.1
$$\begin{array}{l}
\displaystyle n(P_i,P_j)=
\alpha_i+\alpha_j+\sigma_0+\sigma_k\equiv {n-3\choose d-1}\pmod 2\cr
\displaystyle n(P_j,P_k)=\alpha_j+\alpha_k+\sigma_0+\sigma_i\equiv
{n-3\choose d-1}\pmod 2\end{array}$$
(where $\alpha_*$ and $\sigma_*$ are of course as above)
and adding these two equalities we get
$$\alpha_i+\alpha_k+\sigma_i+\sigma_k\equiv 0\pmod 2\ .$$
The identity
$$\sigma_0+\sigma_i+\sigma_j+\sigma_k={n-3\choose d-1}$$
yields now
$$n(P_i,P_k)=
\alpha_i+\alpha_k+\sigma_j+\sigma_0\equiv {n-3\choose d-1}\pmod 2\ .$$
This shows $P_i\sim P_k$ and establishes the transitivity of $\sim$.

In order to prove that the Orchard relation has at most two classes,
consider $P_i\not\sim P_j$ and $P_k\not\sim P_j$. The above
computation shows that $P_i\sim P_k$.
\hfill QED

{\bf Proof of Proposition 1.3.} Consider $P(t),Q(t)\in {\mathcal
  P}(t)$. A hyperplane $H(t)$ determined by $d$ points
  $P_{i_1}(t),\dots,P_{i_d}(t)$ of ${\mathcal
  P}(t)\setminus\{P(t),Q(t)\}$ changes its incidence with
the intervall $[P(t),Q(t)]$ at time $t=0$ if and only if 
$P_{i_1}(t),\dots,P_{i_d}(t)\subset {\mathcal F}(t)$ and exactly
one of the points $P(t),Q(t)$ is the last remaining point of ${\mathcal
  F}(t)$.
This proves the equality
$$n(P(1),Q(1))=n(P(-1),Q(-1))$$
if $P(t),Q(t)\in{\mathcal F}(t)$ or if $P(t),Q(t)\in
{\mathcal P}(t)\setminus {\mathcal F}(t)$. In the remaining case
where $P(t)\in {\mathcal F}(t)$ and 
$Q(t)\in {\mathcal P}(t)\setminus {\mathcal F}(t)$ (up to permutation
of $P$ and $Q$) we have
$$n(P(1),Q(1))=n(P(-1),Q(-1))\pm 1$$
(the difference coming of course from the unique hyperplan spanned by 
${\mathcal F}(t)\setminus P(t)$). Proposition 1.3 follows
now easily.\hfill QED

{\bf Proof of Corollary 1.4.} A 
flip induces an exchange of exactly $(d+1)$ points beetween the two
equivalence classes. This implies assertions (i) and (ii) at once. 
For assertion (iii) we need also the fact that every pair of generic
configurations having $n$ points in ${\bf R}^d$ can be joined by 
a finite number of flips.
\hfill QED

The formula for $\pi(2m,2d+1)$ will follow from Example 3.2 treated in
the next section.

{\bf Proof of Corollary 1.6.} Given a $(d+1)-$dimensional 
vectorspace ${\mathcal C}={\mathcal C}_0+{\bf R}$ of continuous functions
${\bf R}^k\longrightarrow {\bf R}$ and a ${\mathcal C}-$generic
set ${\mathcal P}=\{P_1,\dots,P_n\}$ in ${\bf R}^k$, consider the
{\it generalised Veronese-map} $V:{\bf R}^k\longrightarrow {\bf R}^d$ defined by
$V(P)=(b_1(P),\dots,b_d(P))$ where $b_1,\dots,b_d\in{\mathcal C}_0$
form a basis of the vectorspace ${\mathcal C}_0$. 
This map is of course well-defined up to a linear automorphism
of ${\bf R}^d$ and sends $\mathcal P$ into a
subset $V({\mathcal P})\subset {\bf R}^d$ which is generic in the sense of 
Theorem 1.1: A point $V(P)\in {\bf R}^k$ belongs to a hyperplan $H$
spanned by $V(P_{i_1},\dots,V_{i_d})$ if and only if 
$f(P)=0$ for any function $f\in I(P_{i_1},\dots,P_{i_d})$. 
Such a non-zero function $f$ changes the sign
according  to the two connected components of 
${\bf R}^k\setminus H=H_+\cup H_-$: One has $f(P)<0$ for $P\in{\bf R}^k$ if
and only if $V(P)\in H_-$. Theorem 1.1. implies now 
obviously the result.\hfill QED 

\section{Computational aspects}

Given a finite, totally ordered set $S$, we denote by 
${S\choose k}$ the set of all distinct strictly 
increasing sequences
$$s_{i_1}<s_{i_2}<\dots<s_{i_k}$$
of length $k$ in $S$. 

Order the points ${\mathcal P}=\{P_1,\dots,P_n\}\subset {\bf R}^d$ of
a 
generic configuration
totally (for instance by setting $P_i<P_j$ if $i<j$). For 
$X\in {\bf R}^d$ and $S=\{P_{i_1},\dots,P_{i_d}\}\in {{\mathcal P}
\choose d}$ we define $\hbox{det}(S-X)$ as the determinant of the square 
$d\times d$ matrix with rows
$P_{i_1}-X,\dots,P_{i_d}-X$.

{\bf Proposition 3.1.} {\sl Let ${\mathcal P}=\{P_1,\dots,P_n\}\subset 
{\bf R}^d$
be a generic configuration of $n$ points. We have then
$P\sim Q$ if and only if
$$(-1)^{n-3\choose d-1}\prod_{S\in {{\mathcal P}\setminus \{P,Q\}\choose d}}
\hbox{det}(S-P)\ \hbox{det}(S-Q)\ >0\ .$$
}

{\bf Proof.} The set ${{\mathcal P}\setminus \{P,Q\}\choose d}$
corresponds to the set of all hyperplanes spanned by $d$ points
in ${\mathcal P}\setminus\{P,Q\}$. Such a hyperplan $H=H_S$
separates the points $P$ and $Q$ if and only if 
$\hbox{det}(S-P)\ \hbox{det}(S-Q)<0$.\hfill QED

In particular, the Orchard relation $\sim$ on a generic configuration
$\mathcal P$ of $n$ points in ${\bf R}^d$ can be constructed by computing
$$(n-1){n-2\choose d-1}$$
determinants of $d\times d$ matrices. The
computational cost of determining $\sim$ is hence of order
$O(n^d)$ in $n$ (for fixed $d$).

{\bf Example 3.2.} Choose $0<a_1<a_2<\dots,<a_n$ and set 
$$P_k=(a_k^1,a_k^2,\dots,a_k^d)\in{\bf R}^d$$
for $1\leq k\leq n$. This yields a generic configuration $
\mathcal P$ of $n$ points in ${\bf R}^d$. A simple computation
using Vandermonde's formula shows then that the number
$n(P_k,P_{k+1})$ of hyperplanes separating $P_k$ from $P_{k+1}$
is always zero. The Orchard relation on $\mathcal P$ is hence
trivial if ${n-3\choose d-1}\equiv 0\pmod 2$ and has two 
non-empty classes (for $n\geq 2$) otherwise. This implies easily
the formula for $\pi(2n,2d+1)$ given after Corollary 1.4.

{\bf Remark 3.3.} For practical purposes, the Orchard relation
(for a huge number $n$ of generic points $P_1,\dots,P_n\in {\bf R}^d$) 
is probably best determined as follows: For each affine hyperplane $H$
spanned by $d$ points of $\{P_2,\dots,P_n\}$ compute an affine function 
$f_H$ satisfying $f_H\vert_H\equiv 0$ and $f_H(P_1)=-1$. The hyperplane
$H$ contributes then $1$ to $n(P_1,P_i)$ if and only if $f_H(P_i)>0$
and the knowledge of the parity of all numbers $n(P_1,P_i)$ determines
of course the Orchard relation by transitivity.

The Orchard relation $\sim$ (or the whole binary
rooted tree obtained by recursive iterations of $\sim$ on equivalence
classes) is probably one of the simplest invariants of generic 
configurations. It can of course be combined with other invariants
(e.g. the subset of vertices forming the convex hull) or used to
define other new invariants. 

Such a new invariant is for instance the function
$$\varphi(P)=\hbox{sign}\Big(\prod_{1\leq i_1<i_2<\dots<i_d\leq n,\
  P_{i_j}\not\sim P}\hbox{det}(P_{i_1}-P,\dots,P_{i_d}-P)\Big)\in
\{\pm 1\}$$
on an equivalence class ${\mathcal A}\subset {\mathcal P}=\{P_1,\dots,
P_n\}$.
This function is well-defined if ${n-2-\vert{\mathcal A}\vert\choose
  d-2}
\equiv 0\pmod 2$ and is defined up to a global sign change otherwise.

Similarly, the function
$$\omega(P,Q)=\hbox{sign}\Big(\prod_{1\leq i_1<i_2<\dots<i_{d-1}\leq n,\
  P_{i_j}\not\sim P,Q}\hbox{det}(P-Q,P_{i_1}-Q\dots,P_{i_{d-1}}-Q)\Big)$$
is antisymmetric on an equivalence class ${\mathcal A}\subset
{\mathcal P}$. It is well-defined if ${n-3-\vert{\mathcal
    A}\vert\choose d-3}\equiv 0\pmod 2$ and is defined up to a global
sign otherwise. 

{\bf Remark 3.4.} The formula of Proposition 3.1 suggests perhaps
a homological origin for the Orchard relation. Indeed,
the configuration space of $n$ points in generic position
in ${\bf C}^d$ is connected but (generally) not simply connected.
One gets hence fundamental groups $B_n({\bf C}^d)$ (respectively
$P_n({\bf C}^d)$) by considering loops up to isotopy in this space 
(respectively loops not permuting the points). For $d=1$ we get of
course the braid and the pure braid groups.

The abelianisation of the pure group 
$P_n({\bf C}^d)$ is then isomorphic to ${\bf
  Z}^{n\choose d+1}$. Indeed, each
subset $P_{i_0},\dots,P_{i_d}$ of $d+1$ distinct points in ${\mathcal
P}(t)$ determines a homorphism onto ${\bf Z}$ by considering the 
winding number
$$\frac{1}{2\pi\sqrt{-1}}\int_0^1
\hbox{ln}\big(\hbox{det}(P_{i_1}(t)-P_{i_0}(t),\dots,P_{i_d}(t)-P_{i_0}(t))
\big)$$
along a loop ${\mathcal P}(t)=(P_1(t),\dots,P_n(t)),\ t\in [0,1]$,
in the space of generic configurations. These homorphisms are linarly
independent and the intersection of all their kernels is
the derived group. 

In the case of $B_n({\bf C}^d)$ we get a homomorphism into ${\bf Z}$ by
considering the winding number
$$\frac{1}{2\pi\sqrt{-1}}\int_0^1
\hbox{ln}\big(\prod_{1\leq i_0<i_1<\dots<i_d\leq n} \big(\hbox{det}
(P_{i_1}(t)-P_{i_0}(t),\dots,P_{i_d}(t)-P_{i_0}(t))\big)^2\big)\ .$$ 

\section{Generic configurations in real projective spaces}

A configuration of $n$ points ${\mathcal P}=\{P_1,\dots,P_n\}
\subset {\bf R}P^d$
in real projective space is generic if no subset of $(k+1)\leq 
(d+1)$ points
in $\mathcal P$ is contained in a projective subspace of dimension $<k$. 
The injection $i: {\bf R}^d\longrightarrow
{\bf R}P^d$ constructed by gluing an ${\bf R}P^{d-1}$ along the
boundary of ${\bf R}^d$ yields a surjection (which is generally not
injective) from the set of generic configuration in ${\bf R}^d$ (up
to isomorphism) onto the set of 
generic configurations in ${\bf R}P^d$ (up to the obvious natural
notion of isomorphism obtained by allowing the hyperplane 
at infinity of affine configurations to move).
 
The exact statement of the 
projective counterpart of the Orchard Theorem depends unfortunately
on the parity of the binomial coefficient
${n-2\choose d}$ enumerating all relevant hyperplanes not containing
two given points of a generic
configuration. For the easy case we have:

{\bf Theorem 4.1A.} (Projective Orchard Theorem.)
{\sl If ${n-2\choose d}$ is even then there exists a
  natural partition of the points ${\mathcal
  P}=\{P_1,\dots,P_n\}\subset
P{\bf R}^d$ of a generic projective configuration into two classes.

This partition is given by considering the equivalence classes of the
affine configuration obtained after erasing a generic ${\bf R}P^{d-1}$
not intersecting ${\mathcal P}$.}

Before stating the result if the binomial coefficient ${n-2\choose
  d}$
is odd we need to introduce some notations: Given two points $P\not=
  Q$
of a generic configuration ${\mathcal P}\subset {\bf R}P^d$ we denote
by $L_{P,Q}$ the projective line spanned by them. Denote
by $\alpha$, respectively $\beta={n-2\choose d}-\alpha$, the number of
projective $(d-1)-$dimensional subspaces spanned by $d$ points
of ${\mathcal P}\setminus\{P,Q\}$ which intersect the first, respectively
the second, of the two connected components in $L_{P,Q}\setminus \{P,Q\}$.
Let $\gamma\in\{\alpha,\beta\}$ be the unique integer such that
$$\gamma\equiv {n-3\choose d-1}\pmod 2$$
and denote by $I(P,Q)$ the corresponding connected component of
$L_{P,Q}\setminus \{P,Q\}$. Denote by $\Gamma$ the immerged 
complete graph with vertices $\mathcal P$ and edges $I(P,Q)$ for
$P\not= Q\in{\mathcal P}$. We call a continuous application of
a connected graph $G$ into ${\bf R}P^d$ (for $d\geq 2$)
{\it homologically 
trivial} if the induced group homomorphism 
$i_*: \pi_1(G)\longrightarrow {\bf Z}/2{\bf Z}=\pi_1({\bf R}P^d)$
is trivial.  

{\bf Theorem 4.1B.} {\sl The immersion of the complete graph $\Gamma$ into
${\bf R}P^d$ is homologically trivial.}

Before sketching proofs, let us remark that these projective
versions can be applied (modulo point-hyperplane duality in
projective space) to generic arrangements of hyperplanes in ${\bf
  R}P^d$. 
(A finite set ${\mathcal H}=\{H_1,\dots,H_n\}$ of $n$ distinct
hyperplanes in ${\bf R}P^d$ is {\it generic} if the intersection of
any subset of $k\leq d$ hyperplanes in ${\mathcal H}$ is of
codimension $k$.) In the case of a projective generic line arrangment
${\mathcal L}$
in the projective plane ${\bf R}P^2$, the statement of Theorem A can
visually be seen as follows: Two distinct lines $L_1,L_2\in {\mathcal
  L}$ define two connected components in ${\bf R}P^2\setminus 
\{L_1,L_2\}$ containing all the remaining ${n-2\choose 2}$
intersections of distinct lines in ${\mathcal L}\setminus
\{L_1,L_2\}$. If the binomial coefficient ${n-2\choose 2}$ is even,
the numbers of intersections $L_i\cap L_j,\ L_i,L_j\in{\mathcal
  L}\setminus
\{L_1,L_2\}$ contained in each  connected component of ${\bf
  R}P^2\setminus
\{L_1,L_2\}$ have the same parity. The two lines $L_1$ and $L_2$ are 
equivalent if and only if the above 
parity is given by ${n-3\choose 2-1}\equiv n+1\pmod 2$. An interesting
feature of this construction is the fact that it generalises also 
to generic configurations of pseudolines (cf. \cite{G} for the
definition), see the next section for details. 
A still more general setting for considering 
the (projective) Orchard relation seems to be given by a suitable subset of 
oriented matroids (cf. \cite{OM}), perhaps the set of arrangements
of pseudohyperplanes in (projective or affine) real space
of dimension $d$ which are generic in the sense that $k$ distinct
pseudohyperplanes have an intersection of codimension $k$ for 
$k\leq d$. Corollary 1.6 is a step in this direction.

{\bf Sketch of proof for Theorem 4.1A.} Erase a projective subspace $P{\bf R}^{d-1}$ containing
no point of $\mathcal P\subset P{\bf R}^d$. The equivalence relation
of the resulting affine generic configuration is independent of the
choice of the above subspace.\hfill QED  

{\bf Sketch of proof of Theorem 4.1B.} Suppose by contradiction that
$\Gamma$ contains a loop $\lambda$ which is homologically non-trivial in ${\bf
  R}P^2$. Use isotopies and the Orchard Theorem (after suitable 
affine projections) to reduce the number of intersections of $\lambda$ 
with ${\mathcal P}$. The result follows then by induction.\hfill QED

{\bf Proof of Theorem 1.5.} Consider the $2-$fold cover
${\bf S}^d\longrightarrow {\bf R}P^d$ sending a generic antipodal
spherical configuration
${\mathcal P}=\{\pm P_1,\dots,\pm P_n\}$ into a generic projective
configuration ${\overline {\mathcal P}}=\{{\overline P}_1,\dots,{\overline P}_n\}$ of
${\bf R}P^d$. If ${n-2\choose d}$ is even, lift the equivalence relation on ${\overline
{\mathcal P}}$ onto ${\mathcal P}$ in the obvious way. If ${n-2\choose d}$ is odd,
choose an equatorial circle ${\bf S}^{d-1}\subset {\bf S}^d$ which avoids all
points of ${\mathcal P}$. Up to a sign choice, we can now suppose that the points
$P_1,\dots,P_n$ of ${\bf S}^d$ belong all to the same hemisphere (connected component
of ${\bf S}^d\setminus{\bf S}^{d-1}$) which we identify with the affine space
${\bf R}^{d}$ via the central projection $\pi$ sending ${\bf S}^{d-1}$ 
at infinity. The set
$\{\pi(P_1),\dots,\pi(P_n)\}$ is now a generic set of $n$ points in ${\bf R}^{d}$.
Set now $P_i\sim P_j$ (for $i\not= j$) if and only if $\pi(P_i)\sim\pi(P_j)$ and extend
this by $P_i\not\sim -P_i$. Theorem 4.1B above implies that this defines an equivalence 
relation with 2 classes which is independent of the choice of the 
equatorial circle ${\bf S}^{d-1}\subset{\bf S}^d\setminus {\mathcal P}$. 
\hfill QED  

\section{Simple arrangements of pseudolines in ${\bf R}P^2$}

The aim of this section is to work out some consequences of the
Orchard Theorem for simple arrangements of pseudolines in the projective plane ${\bf R}
P^2$.

A {\it pseudoline} in ${\bf R}P^2$ is a simple, smooth curve isotopic to 
a projective line in ${\bf R}P^2$.
An {\it arrangement} of $n$ pseudolines is a finite set ${\mathcal L}=\{L_1,\dots,
L_n\}$ of $n$ pseudolines intersecting each other transversally exactly once. Such
an arrangement is {\it simple} if no triple intersections occur. If a simple 
arrangement is stretchable (isotopic in the obvious sense to a simple arrangement
of projective lines), then Theorem 4.1A, respectively Theorem 4.1B, 
of the preceeding section imply
that there is a natural equivalence relation on the set of (pseudo)lines if 
${n-2\choose 2}$
is even (i.e. if $n\equiv 2\pmod 4$ or $n\equiv 3\pmod 4$),
respectively that the (pseudo)lines carry a natural orientation 
(up to global reversion of all orientations) if ${n-2\choose 2}$ is odd.

All this remains valid for simple arrangements of pseudolines and we
have of course an analogue of the Flip Proposition. As a consequence, the 
polygones (connected components of ${\bf R}P^2\setminus {\mathcal L}$) 
have extra structures according to the classes
or orientations of their sides. 

Two pseudolines $L_i,L_j$ of a simple pseudoline arrangement 
${\mathcal L}=\{L_1,\dots,L_n\}$ define two open digons by considering the
two connected components of ${\bf R}P^2\setminus \{L_i\cup L_j\}$.
The interior of these digons contain $\alpha$ respectively $\beta$ intersections
$L_s\cap L_t$ of distinct pseudolines $L_s\not=L_t$ in
${\mathcal L}\setminus \{L_i,L_j\}$.

For ${n-2\choose 2}=\alpha+\beta$ even, set $L_i\sim L_j$ if 
$\alpha\equiv{n-3\choose 2}\pmod 2$ (or if $L_j=L_i$).

{\bf Theorem 5.1A.} {\sl For ${n-2\choose 2}$ even, the relation $\sim $ defined
above is an equivalence relation on ${\mathcal L}$ into at most two classes.}

Consider now two oriented pseudolines $L_i^o,L_j^o$. They induce an orientation
on the boundary of exactly one digon $D^o$ in ${\bf R}P^2\setminus\{L_i\cup
L_j\}$ and they have disagreeing orientations on the boundary of the 
remaining digon. We call the orientations of the oriented pseudolines $L_i^o$ and
$L_j^o$ {\it compatible} at $L_i\cap L_j$ if the number $\alpha^o$ of
intersections $L_s\cap L_t$ (for $L_s\not=L_t\in{\mathcal L}\setminus\{L_i,L_j\}$)
contained in the digon $D^o$ satisfies 
$$\alpha^o\equiv {n-3\choose 2}\pmod 2\ .$$
An orientation of all pseudolines in ${\mathcal L}$ is {\it compatible} if
it is compatible at $L_i\cap L_j$ for all $L_i\not=L_j$.

{\bf Theorem 5.1B.} {\sl If ${n-2\choose 2}$ is odd, then there exist exactly
two compatible orientations of all pseudolines in a simple pseudoline 
arrangement ${\mathcal L}=\{L_1,\dots,L_n\}\subset {\bf R}P^2$.
These two compatible orientations induce opposite orientations on all 
pseudolines.} 

We call the partition, respectively compatible orientation, on the set
of pseudolines of a simple pseudoline arrangement $\mathcal L$
the {\it Orchard partition},
respectively {\it Orchard orientation}, on $\mathcal L$. 

{\bf Proof of Theorem 5.1A} Consider three pseudolines $L_i,L_j,L_k$ decomposing
the projective plane into four triangles $\Delta_0,\Delta_i,\Delta_j,\Delta_k$ as
in Figure 5 (here we do not care about orientations). 

\medskip
\centerline{\epsfysize4cm\epsfbox{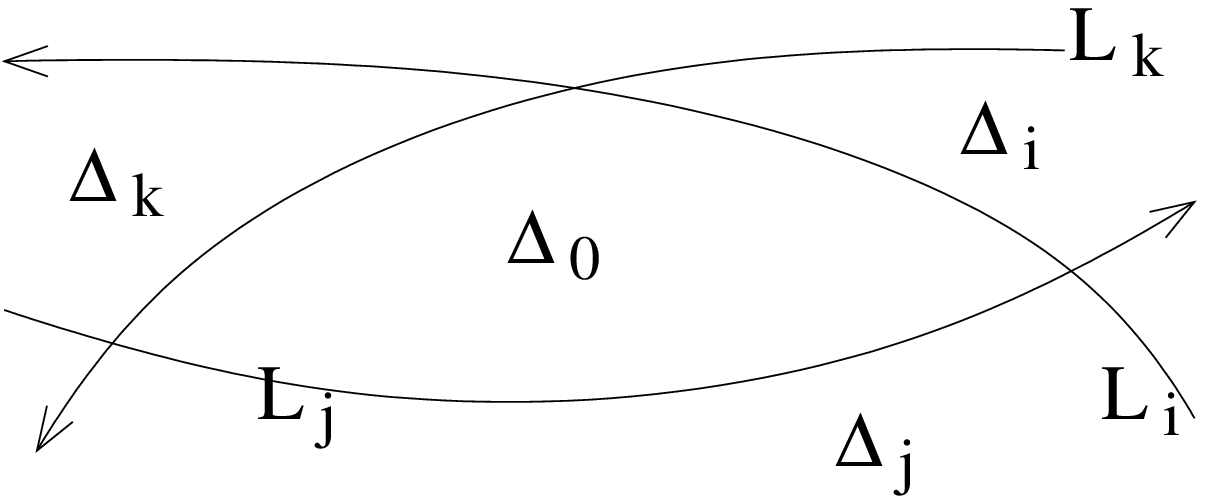}}
\centerline{Figure 5.} 
\medskip

Denote by $\alpha_{i,j}$
the number of pseudolines intersecting $L_i$ and $L_j$ on the
boundary of $\Delta_0$ and define $\alpha_{i,k},\alpha_{j,k}$ similarly.
For $*\in \{0,i,j,k\}$ denote by $\sigma_*$ the number of intersections 
$L_s\cap L_t$ inside $\Delta_*$ of distinct pseudolines $L_s,L_t\in{\mathcal L}
\setminus\{L_i,L_j,L_k\}$. An argument similar to Lemma 2.1 shows that
we have $L_i\sim L_j$ if and only if
$$\sigma_0+\sigma_k+\alpha_{i,k}+\alpha_{j,k}\equiv {n-3\choose 2}\pmod 2$$
and $L_j\sim L_k$ if and only if
$$\sigma_0+\sigma_i+\alpha_{i,j}+\alpha_{i,k}\equiv {n-3\choose 3}\pmod 2\ .$$
The equality 
$$\sigma_0+\sigma_i+\sigma_j+\sigma_k={n-3\choose 3}$$
and arguments similar to those used in the proof of Theorem 1.1
imply now the result.\hfill QED

The proof of Theorem 5.1B is similar and left to the reader.

The analogue of flips are {\it triangle-moves} (corresponding to Reidemeister
III-moves for knots) and they change the Orchard-equivalence relation (respectively
the Orchard-orientation) only on the $3$ pseudolines involved in the
triangle move. Figure 6 shows two simple pseudoline arrangements (with 
arrows indicating one of the two Orchard-orientations) which are 
related by a triangle-move.

\medskip
\centerline{\epsfysize3.5cm\epsfbox{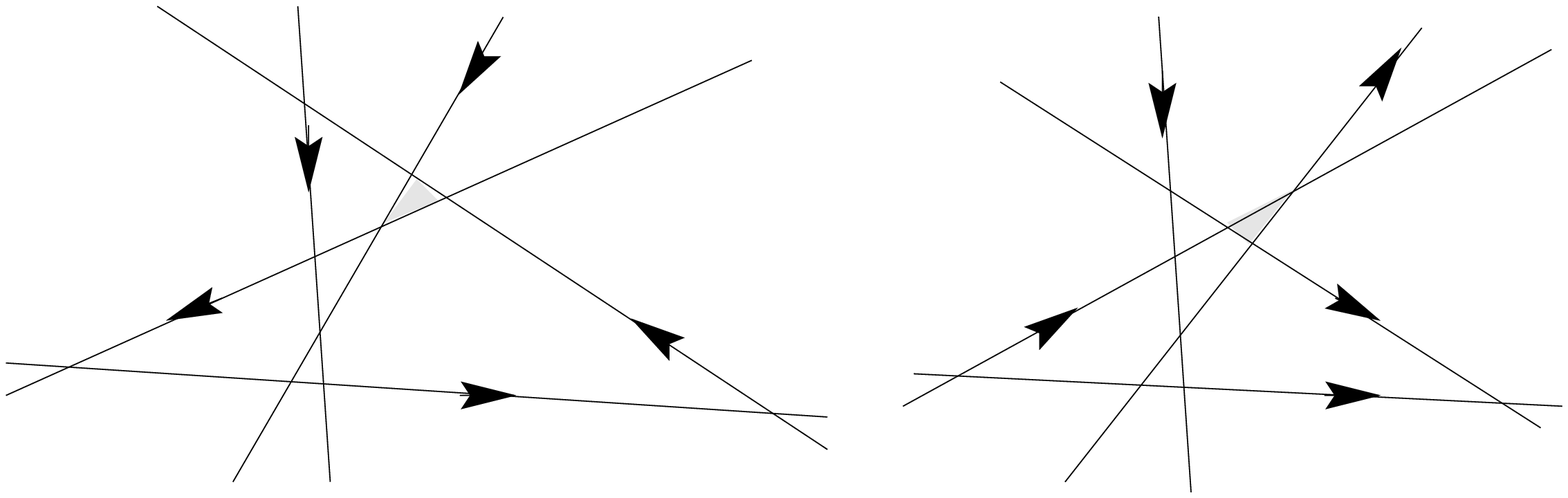}}
\centerline{Figure 6: Two simple (pseudo)line arrangements related by a 
triangle-move.} 
\medskip

Orchard-partitions and Orchard orientations yield of course nice
invariants of pseudoline arrangements. Two such particularly nice invariants
induced by an Orchard orientation (hence in the case where ${n-2\choose 2}$
is odd) are given by desingularising all crossings in one of the two
generic ways (either by respecting all orientations or in the other way)
and by analysing the resulting pattern of noncrossing closed curves
in ${\bf R}P^2$.

I would like to thank many people who where interested in this story and
especially M. Brion, E. Ferrand and A. Marin for their interesting remarks.

\vskip1cm

Roland Bacher

INSTITUT FOURIER

Laboratoire de Math\'ematiques

UMR 5582 (UJF-CNRS)

BP 74

38402 St MARTIN D'H\`ERES Cedex (France)
 
e-mail: Roland.Bacher@ujf-grenoble.fr

\end{document}